\documentclass[10pt]{article}
\usepackage{mathrsfs}
\usepackage{amsfonts}
\usepackage{amsthm,amsmath,amssymb,anysize}
\newtheorem{lemma}{Lemma}[section]
\newtheorem{theorem}[lemma]{Theorem}

\newtheorem{proposition}[lemma]{Proposition}

\usepackage[numbers,sort&compress]{natbib}
\usepackage{amsmath}
\allowdisplaybreaks[3]
\marginsize{4.2cm}{4.2cm}{3.6cm}{3cm}
\numberwithin{equation}{section}

\title{\textsf{Derivations and 2-local derivations on the mirror Heisenberg-Virasoro algebra}}
\author{\textsc{Xuelian Guo and Liming Tang\footnote{corresponding author}}\\
\small{School of Mathematical Sciences}\\
\small{Harbin Normal University}\\
\small{150025 Harbin, China}\\
\small{E-mail: limingtang@hrbnu.edu.cn}}
\date{ }
\date{ }
\begin{document}
\maketitle
 \begin{quotation}
{\small\noindent \textbf{Abstract}:
 Using the first cohomology from the mirror Heisenberg-Virasoro algebra to the twisted Heisenberg algebra (as the mirror Heisenberg-Virasoro algebra-module), in this paper we determined the derivations on the mirror Heisenberg-Virasoro algebra. Based on this result, we proved that any 2-local derivation on the mirror Heisenberg-Virasoro algebra is a derivation.

 \vspace{0.05cm} \noindent{\textbf{Keywords}}: derivation; 2-local derivations; the mirror Heisenberg-Virasoro algebra

\vspace{0.05cm} \noindent \textbf{Mathematics Subject Classification
2010}: 17B05, 17B20, 17B30}
\end{quotation}
 \setcounter{section}{-1}
\section{Introduction}
 Lie algebra is one of the most important tools in physics. For example,  ${\rm Sch\ddot{o}dinger}$-${\rm Virasoro}$ algebras, Heisenberg-Virasoro algebra, twisted Heisenberg-Virasoro algebra and so on are widely used in mathematics and physics. In this paper, we mainly determine  the derivations and 2-local derivations on the Mirror Heisenberg-Virasoro algebra. The mirror  Heisenberg-Virasoro algebra,  whose structure is similar to that of the twisted Heisenberg-Virasoro algebra, is the even part of mirror N = 2 superconformal algebra (see \cite{BK}). As is well known, representation theory is  essential in the research of Lie algebras. The Whittaker modules, $\mathcal{U}(\mathbb{C}d_0)$-free modules  (see (\ref{1.1}), $d_0$ is one of the standard basis vectors for the mirror Heisenberg-Virasoro algebra), and the tensor products of Whittaker modules and $\mathcal{U}(\mathbb{C}d_0)$-free modules on the mirror Heisenerg-Virasoro algebra, which are non-weight modules, are determined. And sufficient and necessary conditions are given for these non-weight modules on the mirror Heisenberg-Virasoro algebra to be irreducibie (see \cite{WE}). Also, the tensor product weight modules on the mirror Heisenberg-Virasoro algebra are studied, and some examples of irreducible weight modules on the mirror Heisenberg-Virasoro algebra are given (see \cite{EH}). All Harish-Chandra modules  on  the mirror Heisenberg-Virasoro algebra are classified (see \cite{ET}). Likewise, structure theory is also very important in the research of Lie algebras. For example, the derivation algebra and automorphism group on the twisted Heisenberg-Virasoro algebra have been determined (see \cite{EM}). 2-local derivations on the twisted  Heisenberg-Virasoro algebra are also given in \cite{E3}.  The mirror  Heisenberg-Virasoro algebra is defined as $\frac{1}{2}\mathbb{Z}$-graded algebra and its $\frac{1}{2}\mathbb{Z}$-graded derivation algebras are given in  \cite{ES}. Different from the graded structure proposed by reference \cite{ES}, in this paper, by virtue of $\mathbb{Z}$-graded structure on the mirror  Heisenberg-Virasoro algebra introduced by \cite{TY}, the first cohomology from the mirror Heisenberg-Virasoro algebra to its ideal the twisted Heisenberg algebra (as the mirror Heisenberg-Virasoro algebra-module) is determined. Then using  the exact sequence of low degree terms associated with  Hochschild-Serre spectral sequence
 \[{\rm E}_{2}^{pq}={\rm H}^{p}(\mathfrak{D}/ \mathcal{H}, {\rm H}^{q}(\mathcal{H}, \mathfrak{D}/\mathcal{H})) \Rightarrow\\[-5pt] {\rm H}^{p+q}(\mathfrak{D}, \mathfrak{D}/ \mathcal{H}),\]
 the exact sequence related to $\mathcal{H}$
 \[0\rightarrow \mathcal{H} \rightarrow \mathfrak{D} \rightarrow \mathfrak{D}/ \mathcal{H} \rightarrow 0\]
 and the exact sequence
 \[{\rm H}^{1}(\mathfrak{D}, \mathcal{H})\rightarrow {\rm H}^{1}(\mathfrak{D}, \mathfrak{D})\rightarrow {\rm H}^{1}(\mathfrak{D}, \mathfrak{D}/\mathcal{H}),\]
  the first cohomology from the mirror Heisenberg-Virasoro algebra to itself is equal to the first cohomology from the mirror Heisenberg-Virasoro algebra to its ideal the twisted Heisenberg algebra. Hence, the $\mathbb{Z}$-graded derivations on the mirror Heisenberg-Virasoro algebra are determined.  Furthermore, based on the result, we proved that every 2-local derivations on the mirror Heisenberg-Virasoro algebra is a derivation.


\section{Derivations on the Mirror Heisenberg-Virasoro algebra}

Hereafter, we denote by $\mathbb{C}$, $\mathbb{Z}$ the set of all complex numbers, integers, respectively. All vector spaces, algebras and their tensor products are assumed to be over  the field $\mathbb{C}$.  Kronecker delta $\delta_{i,j}$ denote 1 if $i=j$ and 0 otherwise.
 We consider the following Lie algebras, which are referred to as $\emph{the Mirror Heisenberg-Virasoro algebra}$ $\frak{D}$  studied by  \cite{WE, EH, ET, TY, BK} . The Lie algebra has $\mathbb{C}$-basis
 \begin{align}\label{1.1}
 &\{d_{m}, h_{r}, {\bf c, l}\ |\ m\in \mathbb{Z}, r\in \frac{1}{2}+ \mathbb{Z}\}
 \end{align}
 with the following nontrivial Lie brackets:
\begin{align*}
&[d_{m}, d_{n}]=(m-n)d_{m+n}+\frac{m^{3}-m}{12}\delta_{m+n}, _{0}{\bf c},\\
&[d_{m}, h_{r}]=-rh_{m+r},\\
&[h_{r}, h_{s}]=r\delta_{r+s,0}{\bf l},\\
&[{\bf c}, \mathfrak{D}]=[{\bf l}, \mathfrak{D}]=0,
\end{align*}
where  $m, n\in \mathbb{Z}$, $r,s\in \frac{1}{2}+\mathbb{Z}.$ From \cite{TY}, suppose ${\rm deg}d_{n}=n$, ${\rm deg}h_{r}=r-\frac{1}{2}$, ${\rm deg}{\bf c}=0$, ${\rm deg}{\bf l}=-1$.
  Let
\begin{align*}
&\mathfrak{D}_{0}={\rm Span}_{\mathbb{C}}\{d_{0},h_{\frac{1}{2}},{\bf c}\},\\
&\mathfrak{D}_{-1}={\rm Span}_{\mathbb{C}}\{d_{-1},h_{-\frac{1}{2}}, {\bf l}\},\\
&\mathfrak{D}_{n}={\rm Span}_{\mathbb{C}}\{d_{n},h_{n+\frac{1}{2}}\}, n\in \mathbb{Z}\backslash\{0,-1\}.
\end{align*}
Then $\mathfrak{D}=\oplus_{{n\in \mathbb{Z}}}\mathfrak{D}_{n}$ is a ${\mathbb{Z}}$-graded algebra.
 Note that $\frak{D}$ is a generalization of the $\emph{Virasoro}$ algebra $\mathcal{V}$ and $\emph{the twisted Heisenberg}$ algebra $\mathcal{H}$. Furthermore, the Lie algebra spanned by $\{d_{m}, {\bf c}\ |\ m\in \mathbb{Z}\}$ is isomorphic to the $\emph{Virasoro}$ algebra $\mathcal{V}$ and  the Lie algebra spanned by $\{h_{r}, {\bf l}\ |\ r\in \frac{1}{2}+ \mathbb{Z}\}$ is isomorphic to $\emph{the twisted Heisenberg}$ algebra $\mathcal{H}$. It is necessary to remind that $\mathcal{H}$ is the ideal of $\mathfrak{D}$ and $\mathfrak{D}$ is the semi-direct product of $\mathcal{V}$ and $\mathcal{H}$, that is,\ $\mathfrak{D}=\mathcal{V}\ltimes\mathcal{H}$. The center of $\mathfrak{D}$  is spanned by $\{\bf c,l\}$, denoted by ${\rm{C}}(\mathfrak{D})$. Let $\mathcal{U}(\mathcal{V})$ denote the universal enveloping algebra of $\mathcal{V}$.

 \begin{proposition}
   Let $\mathfrak{D}=\oplus_{n\in\mathbb{Z}}\mathfrak{D}_n$ be $\emph{the mirror Heisenberg-Virasoro}$ algebra and  $\mathcal{H}$   $\emph{the twisted Heisenberg}$ algebra. Suppose  $\mathcal{H}_{-1}=\mathbb{C}h_{-\frac{1}{2}}\oplus \mathbb{C}{\bf l}$ and  $\mathcal{H}_{n}=\mathbb{C}h_{n+\frac{1}{2}}$, $n\in \mathbb{Z}\setminus\{-1\}$, then
 \begin{itemize}
   \item[\rm{(i)}] $\mathcal{H}=\oplus_{n\in \mathbb{Z}}\mathbb{C}h_{n+\frac{1}{2}}\oplus \mathbb{C}{\bf l}$ is a $\mathbb{Z}$-graded Lie algebra.
   \item[\rm{(ii)}] $\mathcal{H}$ is a $\mathbb{Z}$-graded $\mathfrak{D}$-module.
 \end{itemize}
 \end{proposition}

  Let $L$ be a Lie algebra and $V$ be an  $L$-module. A linear map $\varphi$  from $L$ to $V$ is called a ${\emph{derivation}}$ if it satisfied $\varphi([x, y])=x\cdot\varphi(y)-y\cdot\varphi(x)$ for any $x, y\in L$. A linear map $\varphi: x\mapsto x\cdot v$ is called an  ${\emph{inner derivation}}$ for any $x\in L$ and $v\in V$.  If $L$ and $V$ are $\mathbb{Z}$-graded and $\varphi(L_{n})\subset V_{m+n}$ for all $m,n\in \mathbb{Z},$  then $\emph{order}$ of $\varphi$ is $m$, denoted by  $\deg\varphi=m$. Denoted by ${\rm Der}_{\mathbb{C}}(L, V)$ the vector space consisting of all derivations from $L$ to $V$. Denoted by ${\rm Inn}_{\mathbb{C}}(L, V)$ the vector space consisting of all inner derivations from $L$ to $V$.

  Let $${\rm H}^{1}(\mathfrak{D}, \mathcal{H})={\rm Der}_{\mathbb{C}}(\mathfrak{D}, \mathcal{H})/{\rm Inn}_{\mathbb{C}}(\mathfrak{D}, \mathcal{H}).$$
By virtue of \cite{CM}, the exact sequence given by $\mathcal{H}$
\[
0\rightarrow \mathcal{H} \rightarrow \mathfrak{D} \rightarrow \mathfrak{D}/ \mathcal{H} \rightarrow 0
\]
induces an exact sequence
\begin{equation}
{\rm H}^{1}(\mathfrak{D}, \mathcal{H})\rightarrow {\rm H}^{1}(\mathfrak{D}, \mathfrak{D})\rightarrow {\rm H}^{1}(\mathfrak{D}, \mathfrak{D}/\mathcal{H}) \label{9}
\end{equation}
of $\mathbb{Z}$-graded vector spaces.
According to \cite{WA}, the right side can be calculated from the exact sequence of low degree terms
\begin{equation}
0 \rightarrow {\rm H}^{1}(\mathfrak{D}/ \mathcal{H}, \mathfrak{D}/ \mathcal{H})\rightarrow {\rm H}^{1}(\mathfrak{D}, \mathfrak{D}/\mathcal{H})\rightarrow {\rm H}^{1}(\mathcal{H}, \mathfrak{D}/\mathcal{H})^{\mathfrak{D}/\mathcal{H}}\label{2}
\end{equation}
associated with Hochschild-Serre spectral sequence  $${\rm E}_{2}^{pq}={\rm H}^{p}(\mathfrak{D}/ \mathcal{H}, {\rm H}^{q}(\mathcal{H}, \mathfrak{D}/\hslash)) \Rightarrow\\[-5pt] {\rm H}^{p+q}(\mathfrak{D}, \mathfrak{D}/ \mathcal{H}).$$
Here, according to
 ${\rm H}^{1}(\mathfrak{D}/ \mathcal{H}, \mathfrak{D}/ \mathcal{H})=0$ (see \cite{LZ}) and
 \[
 {\rm H}^{1}(\mathcal{H}, \mathfrak{D}/\mathcal{H})^{\mathfrak{D}/\mathcal{H}}
 \cong {\rm Hom}_{\mathbb{C}}(\mathcal{H}/[\mathcal{H}, \mathcal{H}], \mathfrak{D}/\mathcal{H})^{\mathfrak{D}/\mathcal{H}}
 ={\rm Hom}_{\mathfrak{D}/\mathcal{H}}(\mathcal{H}/[\mathcal{H}, \mathcal{H}], \mathfrak{D}/\mathcal{H})
 \] (see \cite{WA}),
${\rm H}^{1}(\mathcal{H}, \mathfrak{D}/\mathcal{H})^{\mathfrak{D}/\mathcal{H}}$ can be embedded in ${\rm Hom}_{\mathcal{U}(\mathcal{V})}(\mathcal{H}/[\mathcal{H}, \mathcal{H}], \mathcal{V})$.
\begin{proposition}\label{YY}
 Let $\mathfrak{D}=\oplus_{n\in\mathbb{Z}}\mathfrak{D}_n$ be the mirror Heisenberg-Virasoro algebra and $\mathcal{H}$  the twisted Heisenberg algebra. If $\mathcal{H}=\oplus_{n\in \mathbb{Z}}\mathbb{C}h_{n+\frac{1}{2}}\oplus \mathbb{C}{\bf l}$ is a $\mathbb{Z}$-graded $\mathfrak{D}$-module, then $${\rm Der_{\mathbb{C}}}(\mathfrak{D}, \mathcal{H})=\oplus_{n\in \mathbb{Z}}{\rm Der_{\mathbb{C}}}(\mathfrak{D}, \mathcal{H})_{n},$$ where $${\rm Der_{\mathbb{C}}}(\mathfrak{D}, \mathcal{H})_{n}\colon\!\!\!=\{\varphi\in {\rm Der_{\mathbb{C}}}(\mathfrak{D}, \mathcal{H})\ |\ {\rm deg}\varphi=n\}.$$
\end{proposition}

\begin{lemma}\label{UU}
  Let $\mathfrak{D}=\oplus_{n\in\mathbb{Z}}\mathfrak{D}_n$ be the mirror Heisenberg-Virasoro algebra and $\mathcal{H}=\oplus_{n\in \mathbb{Z}}\mathbb{C}h_{n+\frac{1}{2}}\oplus \mathbb{C}{\bf l}$  the twisted Heisenberg algebra. Then $${\rm H}^{1}(\mathfrak{D}_{0}, \mathcal{H}_{n})=0, \ n\in \mathbb{Z}\setminus\!\{0,-1\}.$$
\begin{proof}
  Since $${\rm H}^{1}(\mathfrak{D}_{0}, \mathcal{H}_{n})={\rm Der}_{\mathbb{C}}(\mathfrak{D}_{0}, \mathcal{H}_{n})/ {\rm Inn}_{\mathbb{C}}(\mathfrak{D}_{0}, \mathcal{H}_{n})$$ for $\varphi\in
  {\rm Der}_{\mathbb{C}}(\mathfrak{D}_{0}, \mathcal{H}_{n})$, then it is enough to verify that $\varphi\in
  {\rm Inn}_{\mathbb{C}}(\mathfrak{D}_{0},\mathcal{H}_{n})$.\\[-5pt]

  For \ $\varphi\in {\rm Der}_{\mathbb{C}}(\mathfrak{D}_{0}, \mathcal{H}_{n})$,
$$\varphi\colon\mathfrak{D}_{0}\rightarrow \mathcal{H}_{n}$$
given by
\begin{align*}
\varphi(d_{0})&=a_{n}h_{n+\frac{1}{2}},\\
\varphi(h_{\frac{1}{2}})&=b_{n}h_{n+\frac{1}{2}},\\
\varphi({\bf c})&=c_{n}h_{n+\frac{1}{2}},
\end{align*}
then
\begin{align*}
\varphi([d_{0}, h_{\frac{1}{2}}])&=[\varphi(d_{0}), h_{\frac{1}{2}}]+[d_{0}
, \varphi(h_{\frac{1}{2}})]\\&=[a_{n}h_{n+\frac{1}{2}}, h_{\frac{1}{2}}]+[d_{0}, b_{n}h_{n+\frac{1}{2}}]\\&=-(n+\frac{1}{2})b_{n}h_{n+\frac{1}{2}}\\&=-\frac{1}{2}\varphi(h_{\frac{1}{2}}).
\end{align*}
Hence \[ \varphi(h_{\frac{1}{2}})=2(n+\frac{1}{2})b_{n}h_{n+\frac{1}{2}}=b_{n}h_{n+\frac{1}{2}}.\]
This concludes $b_{n}=0$, that is $\varphi(h_{\frac{1}{2}})=0$.
In additional,
\begin{align*}
\varphi([d_{0},\ {\bf c}])&=[\varphi(d_{0}),\ {\bf c}]+[d_{0}, \varphi({\bf c})]\\&=[d_{0}, c_{n}h_{n+\frac{1}{2}}]\\&=-c_{n}(n+\frac{1}{2})h_{n+\frac{1}{2}}\\&=0.
\end{align*}
It follows that $c_{n}=0$, that is, $\varphi({\bf c})=0$.\\[-5pt]

Let $E_{n}=-\frac{a_{n}}{n+\frac{1}{2}}h_{n+\frac{1}{2}}$. Then
\begin{align*}
\varphi(d_{0})&=[d_{0}, E_{n}]\\&=[d_{0}, -\frac{a_{n}}{n+\frac{1}{2}}h_{n+\frac{1}{2}}]\\&=a_{n}h_{n+\frac{1}{2}}.
\end{align*}
Thus $\varphi\colon\mathfrak{D}_{0}\rightarrow\mathcal{H}_{n}$ can be given by
\begin{align*}
d_{0}&\mapsto[d_{0}, E_{n}]\\h_{\frac{1}{2}}&\mapsto[h_{\frac{1}{2}}, E_{n}]\\{\bf c}&\mapsto[{\bf c}, E_{n}].
\end{align*}
Hence, ${\rm H}^{1}(\mathfrak{D}_{0}, \mathcal{H}_{n})=0,  n\in \mathbb{Z}\setminus\{0,-1\}$.
\end{proof}
\end{lemma}

\begin{lemma}\label{RR}
   Let $\mathfrak{D}=\oplus_{n\in\mathbb{Z}}\mathfrak{D}_{n}$ be the mirror Heisenberg-Virasoro algebra and $\mathcal{H}=\oplus_{n\in \mathbb{Z}}\mathbb{C}h_{n+\frac{1}{2}}\oplus \mathbb{C}{\bf l}$ the twisted Heisenberg algebra. Then $${\rm Hom}_{\mathfrak{D}_{0}}(\mathfrak{D}_{m}, \mathcal{H}_{n})=0\, \mbox{for all}  \,\, m\neq n, m\neq0.$$
\end{lemma}
\begin{proof}
   Recall that $${\rm Hom}_{\mathfrak{D}_{0}}(\mathfrak{D}_{m}, \mathcal{H}_{n})=\{f\colon\mathfrak{D}_{m}\rightarrow \mathcal{H}_{n}\ |  x\cdot f(y)=f(x\cdot y), x\in \mathfrak{D}_{0}, y\in \mathfrak{D}_{m}\}.$$
  \item[Case 1:] If $m\in \mathbb{Z}\backslash\{0, -1\}$, then $$\mathfrak{D}_{m}={\rm Span}_{\mathbb{C}}\{d_{m}, h_{m+\frac{1}{2}}\}.$$

When $n\neq -1$, supoose $f(d_{m})=kh_{n+\frac{1}{2}},\ k\in \mathbb{C},$ then
\begin{align*}
[d_{0}, f(d_{m})]&=f([d_{0}, d_{m}])\\&=-mf(d_{m})\\&=-mkh_{n+\frac{1}{2}}\\&=[d_{0}, kh_{n+\frac{1}{2}}]\\&=-k(n+\frac{1}{2})h_{n+\frac{1}{2}}.
\end{align*}
Hence $k=0.$

Suppose $f(h_{m+\frac{1}{2}})=kh_{n+\frac{1}{2}},\ k\in \mathbb{C},$ then
\begin{align*}
[d_{0}, f(h_{m+\frac{1}{2}})]&=f([d_{0}, h_{m+\frac{1}{2}}])
\\&=-(m+\frac{1}{2})f(h_{m+\frac{1}{2}})
\\&=-(m+\frac{1}{2})kh_{n+\frac{1}{2}}
\\&=-(n+\frac{1}{2})kh_{n+\frac{1}{2}}.
\end{align*}
Hence $k=0$.

When $n=-1$, suppose $$f(d_{m})=k_{1}h_{-\frac{1}{2}}+k_{2}{\bf l},\  f(h_{m+\frac{1}{2}})=k_{3}h_{-\frac{1}{2}}+k_{4}{\bf l},\ k_{i}\in \mathbb{C}, i=1,2,3,4,$$ then
\begin{align*}
[d_{0}, f(d_{m})]&=f([d_{0}, d_{m}])\\&=[d_{0}, k_{1}h_{-\frac{1}{2}}+k_{2}{\bf l}]\\&=\frac{1}{2}k_{1}h_{-\frac{1}{2}}\\&=-mk_{1}h_{-\frac{1}{2}}-mk_{2}{\bf l}.
\end{align*}
Hence $k_{1}=k_{2}=0.$
Furthermore, we have
\begin{align*}
[d_{0}, f(h_{m+\frac{1}{2}})]&=f([d_{0}, h_{m+\frac{1}{2}}])\\&=[d_{0}, k_{3}h_{-\frac{1}{2}}+k_{4}{\bf l}]\\&=\frac{1}{2}k_{3}h_{-\frac{1}{2}}\\&=-(m+\frac{1}{2})k_{3}h_{-\frac{1}{2}}-(m+\frac{1}{2})k_{4}{\bf l}.
\end{align*}
Hence $k_{3}=k_{4}=0.$
\item[Case 2:] When $m=-1$, since $m\neq n$, then $n\neq-1$. Here $$\mathfrak{D}_{-1}={\rm Span}_{\mathbb{C}}\{d_{-1}, h_{-\frac{1}{2}},{\bf l}\}.$$
Let $$f(d_{-1})=k_1h_{n+\frac{1}{2}},\  f(h_{-\frac{1}{2}})=k_2h_{n+\frac{1}{2}},\ f({\bf l})=k_3h_{n+\frac{1}{2}}\,\,\mbox{for}\,  k_1, k_2, k_3\in \mathbb{C}.$$ Then we have
$$[d_{0}, f(d_{-1})]=f([d_{0}, d_{-1}])=f(d_{-1})=k_1h_{n+\frac{1}{2}}=k_1[d_{0},h_{n+\frac{1}{2}}]=-k_1(n+\frac{1}{2})h_{n+\frac{1}{2}}.$$
So $k_1=0$. Furthermore, we have
$$[d_{0}, f(h_{-\frac{1}{2}})]=f([d_{0}, h_{-\frac{1}{2}}])=\frac{1}{2}k_2h_{n+\frac{1}{2}}=k_2[d_{0}, h_{n+\frac{1}{2}}]=-(n+\frac{1}{2})k_2h_{n+\frac{1}{2}}.$$
Hence $k_2=0$. In additional, we have
$$[d_{0},f({\bf l})]=f([d_{0}, {\bf l}])=k_3[d_{0}, h_{n+\frac{1}{2}}]=-k_3(n+\frac{1}{2})h_{n+\frac{1}{2}}=0.$$
Hence $k_3=0$.

Thus $f=0$. That is ${\rm Hom}_{\mathfrak{D}_{0}}(\mathfrak{D}_{m}, \mathcal{H}_{n})=0$ for $m\neq n, m\neq0$.
\end{proof}

\begin{lemma}\label{PP}
   Let $\mathfrak{D}=\oplus_{n\in\mathbb{Z}}\mathfrak{D}_{n}$ be the mirror Heisenberg-Virasoro algebra and $\mathcal{H}=\oplus_{n\in \mathbb{Z}}\mathbb{C}h_{n+\frac{1}{2}}\oplus \mathbb{C}{\bf l}$\ the twisted Heisenberg algebra.  If $m=0$ and $n\neq0,-1$, then ${\rm Hom}_{\mathfrak{D}_{0}}(\mathfrak{D}_{0}, \mathcal{H}_{n})=0$.
\begin{proof}
 Suppose $f(d_{0})=k_1h_{n+\frac{1}{2}}, f(h_\frac{1}{2})=k_2h_{n+\frac{1}{2}}, f(c)=k_3h_{n+\frac{1}{2}}$, then
\begin{align*}
f([d_{0}, d_{0}])&=[d_{0}, f(d_{0})]\\&=k_1[d_{0}, h_{n+\frac{1}{2}}]\\&=-k_1(n+\frac{1}{2})h_{n+\frac{1}{2}}\\&=0.
\end{align*}
Hence $k_1=0$. Furthermore, we have
\begin{align*}
f([d_{0}, h_{\frac{1}{2}}])&=[d_{0}, f(h_\frac{1}{2})]
\\&=k_2[d_{0},h_{n+\frac{1}{2}}]
\\&=-k_2(n+\frac{1}{2})h_{n+\frac{1}{2}}
\\&=-\frac{1}{2}f(h_{\frac{1}{2}})
\\&=-\frac{1}{2}k_2h_{n+\frac{1}{2}}.
\end{align*}
Hence $k_2=0$. In additional, we have
\begin{align*}
f([d_{0}, c])=[d_{0}, f(c)]=k_3[d_{0}, h_{n+\frac{1}{2}}]=-k_3(n+\frac{1}{2})h_{n+\frac{1}{2}}
\end{align*}
Hence $k_3=0$.\\[-5pt]

 So $f=0$, that is,  ${\rm Hom}_{\mathfrak{D}_{0}}(\mathfrak{D}_{0}, \mathcal{H}_{n})=0$ for $n\neq0,-1$.
\end{proof}
\end{lemma}
\begin{lemma}\label{2.10}
 Let $\mathfrak{D}$ be the mirror Heisenberg-Virasoro algebra and $\mathcal{H}$\ the twisted Heisenberg algebra. Then $D\in{\rm Der}_{\mathbb{C}}(\mathfrak{D}, \mathcal{H})_{-1}$ if and only if $D$ is a linear map from $\mathfrak{D}$ to $\mathcal{H}$ satisfying the following conditions:
\begin{itemize}
  \item[{\rm{(i)}}] $D(d_{n})=ch_{n-\frac{1}{2}},$
  \item[{\rm{(ii)}}] $D({\bf c})=0,$
  \item[{\rm{(iii)}}] $D(h_{n+\frac{1}{2}})=\delta_{n,0}c{\bf l},$
  \item[{\rm{(iv)}}] $D({\bf l})=0$,
\end{itemize}
where $n\in\mathbb{Z},\ c\in\mathbb{C}.$ Furthermore, D is a inner derivation.
\end{lemma}
\begin{proof}
Sufficiency: Let $D$ be a linear map that satisfies (i)-(iv). Then it is easy to prove that $D$ is a derivation from $\mathfrak{D}$ to $\mathcal{H}$\\[-5pt]

Necessity:
  (i) For $n\neq 0$, set $D(d_{n})=f(n)h_{n-\frac{1}{2}}$. If $n\neq \pm m$, on the one hand
\begin{align*}
  D([d_{m}, d_{n}])&=(m-n)D(d_{m+n})
                 \\&=(m-n)f(m+n)h_{m+n-\frac{1}{2}},
\end{align*}
on the other hand
\begin{align*}
D([d_{m}, d_{n}])&=[D(d_{m}), d_{n}]+[d_{m}, D(d_{n})]
\\[5pt]&=f(m)[h_{m-\frac{1}{2}}, d_{n}]+f(n)[d_{m}, h_{n-\frac{1}{2}}]
\\&=(m-\frac{1}{2})f(m)h_{m+n-\frac{1}{2}}-(n-\frac{1}{2})f(m)h_{m+n-\frac{1}{2}},
\end{align*}
then
\[
(m-n)f(m+n)=(m-\frac{1}{2})f(m)-(n-\frac{1}{2})f(n).
\]
Let $f(1)=f(-1)=c, c\in \mathbb{C}$. Using induction on $n$, assuming $n-1$, the conclusion holds, that is $f(n-1)=c$. Then
\begin{align*}
  (n-2)f(n)&=(n-\frac{3}{2})f(n-1)-\frac{1}{2}f(1)
         \\&=(n-2)c.
\end{align*}
Therefore, $f(n)=c$.
Since $$D([d_{1}, d_{-1}])=2D(d_{0})$$
and
\begin{align*}
D([d_{1}, d_{-1}])&=[D(d_{1}), d_{-1}]+[d_{1}, D(d_{-1})]
\\&=c[h_{\frac{1}{2}}, d_{-1}]+c[d_{1}, h_{-\frac{3}{2}}]
\\&=2ch_{-\frac{1}{2}}.
\end{align*}
Hence, $D(d_{0})=ch_{-\frac{1}{2}}$. So we have $D(d_{n})=ch_{n-\frac{1}{2}}$ for $n\in \mathbb{Z}$.

(ii) For $n\neq0$, let $D(h_{n+\frac{1}{2}})=g(n)h_{n-\frac{1}{2}}$.
If $m\neq -n, n\neq 0$, on the one hand
\begin{align*}
  D([d_{m}, h_{n+\frac{1}{2}}])&=-(n+\frac{1}{2})D(h_{n+m+\frac{1}{2}})
                           \\&=-(n+\frac{1}{2})g(n+m)h_{n+m-\frac{1}{2}},
\end{align*}
on the other  hand
\begin{align*}
D([d_{m}, h_{n+\frac{1}{2}}])&=[D(d_{m}), h_{n+\frac{1}{2}}]+[d_{m}, D(h_{n+\frac{1}{2}})]
                         \\[5pt]&=c[h_{m-\frac{1}{2}}, h_{n+\frac{1}{2}}]+g(n)[d_{m}, h_{n-\frac{1}{2}}]
                         \\&=-(n-\frac{1}{2})g(n)h_{n+m-\frac{1}{2}}.
\end{align*}
Therefore, $(n+\frac{1}{2})g(n+m)=(n-\frac{1}{2})g(n)$.

Let $g(1)=d,\ d\in \mathbb{C}$. Then $g(n)=\frac{d}{2n-1}.$ Using induction on $n$, if $n-1$, the conclusion holds, that is $g(n-1)=\frac{d}{2(n-1)-1}$. Then
\begin{align*}
(n-\frac{1}{2})g(n)&=(n-\frac{3}{2})g(n-1)
                 \\&=(n-\frac{3}{2})\frac{d}{2n-3}.
\end{align*}
Hence $g(n)=\frac{d}{2n-1}$.
And since  $$D([d_{1}, h_{-\frac{1}{2}}])=\frac{1}{2}D(h_\frac{1}{2})$$
and
\begin{align*}
D([d_{1}, h_{-\frac{1}{2}}])&=[D(d_{1}), h_{-\frac{1}{2}}]+[d_{1}, D(h_{-\frac{1}{2}})]
\\&=c[h_\frac{1}{2}, h_{-\frac{1}{2}}]+g(-1)[d_{1}, h_{-\frac{3}{2}}]
\\&=\frac{1}{2}c{\bf l}-\frac{1}{2}dh_{-\frac{1}{2}}.
\end{align*}
Then $D(h_{\frac{1}{2}})=c{\bf l}-dh_{-\frac{1}{2}}$. So we have $$D(h_{n+\frac{1}{2}})=\frac{d}{2n-1}h_{n-\frac{1}{2}}+\delta_{n},_{0}c{\bf l}.$$
Again, on the one hand
\begin{align*}
  D([d_0, h_{\frac12}])&=[D(d_0), h_\frac12]+[d_0, D(h_\frac12)]\\
                       &=[ch_{-\frac12}, h_\frac12]+[d_0, c{\bf l}-dh_{-\frac12}]\\
                       &=-\frac12c{\bf l}-\frac12dh_{-\frac12}
\end{align*}
on the other hand
\[ D([d_0, h_{\frac12}])=-\frac12D(h_{\frac12})=-\frac12c{\bf l}+\frac12dh_{-\frac12}.\]
Comparing the two formulas above, we have $d=0.$ Therefore, $D(h_{n+\frac{1}{2}})=\delta_{n,0}c{\bf l}.$

(iii) Since
\begin{align*}
D([d_{2}, d_{-2}])&=4D(d_{0})+\frac{1}{2}D({\bf c})
                 \\&=4ch_{-\frac{1}{2}}+\frac{1}{2}D({\bf c})
\end{align*}
and
\begin{align*}
D([d_{2}, d_{-2}])&=[D(d_{2}), d_{-2}]+[d_{2}, D(d_{-2})]
\\[5pt]&=c[h_{\frac{3}{2}}, d_{-2}]+c[d_{2}, h_{-\frac{5}{2}}]
\\&=4ch_{-\frac{1}{2}},
\end{align*}
then we have $D({\bf c})=0$.

 (iv) On the one hand  $$D([h_{\frac{1}{2}}, h_{-\frac{1}{2}}])=\frac{1}{2}D({\bf l}),$$
on the other hand
\begin{align*}
D([h_{\frac{1}{2}}, h_{-\frac{1}{2}}])&=[D(h_{\frac{1}{2}}), h_{-\frac{1}{2}}]+[h_\frac{1}{2}, D(h_{-\frac{1}{2}})]
\\&=[c{\bf l}, h_{-\frac{1}{2}}]
\\&=0.
\end{align*}
Then we have $D({\bf l})=0$.

Let $D$ be given by
\begin{align*}
&D(d_n)=[d_n, 2ch_{-\frac12}]=ch_{n-\frac12},\\
&D(h_{n+\frac12})=[h_{n+\frac12}, 2ch_{-\frac12}]=c\delta_{n, 0}{\bf l},\\
&D({\bf c})=[{\bf c}, 2ch_{-\frac12}]=0,\\
&D({\bf l})=[{\bf l}, 2ch_{-\frac12}]=0.
\end{align*}
Then $D$\ is a inner derivation.
\end{proof}

\begin{theorem}\label{QQ}
 Let $\mathfrak{D}=\oplus_{n\in\mathbb{Z}}\mathfrak{D}_{n}$ be the mirror Heisenberg-Virasoro algebra and $\mathcal{H}$\ the twisted Heisenberg algebra. Then
$${\rm Der}_{\mathbb{C}}(\mathfrak{D}, \mathcal{H})={\rm Der}_{\mathbb{C}}(\mathfrak{D}, \mathcal{H})_{0}\oplus{\rm Inn}_{\mathbb{C}}(\mathfrak{D}, \mathcal{H}).$$
\end{theorem}
\begin{proof}
Suppose $$\varphi\colon\mathfrak{D}\rightarrow \mathcal{H}$$ is a derivation, according to proposition \ref{YY},\ then $\varphi$\ can be decomposed into homogeneous elements: \[
\varphi=\Sigma_{n\in \mathbb{Z}}\varphi_{n}\]
where
$$\varphi_{n}\in {\rm Der}_{\mathbb{C}}(\mathfrak{D}, \mathcal{H})_{n}.$$

Let $n\in \mathbb{Z}\setminus\{0,-1\}$. Then $\varphi_{n}|_{\mathfrak{D}_{0}}$ is a derivation of $\mathcal{H}_{n}$  from $\mathfrak{D}_{0}$ to $\mathfrak{D}_{0}$-module. By virtue of Lemma \ref{UU}, $\varphi_{n}|_{\mathfrak{D}_{0}}$ is inner derivation, that is, there exists $E_{n}\in \mathcal{H}_{n}$ such that $\varphi_{n}(x)=[x, E_{n}],\ x\in \mathfrak{D}_{0}$. Considering $\psi_{n}(x)=\varphi_{n}(x)-[x, E_{n}], x\in \mathfrak{D}$, then $\psi_{n}$\ is an inner derivation, ${\rm deg}\psi_{n}=n$, and $\psi_{n}(\mathfrak{D}_{0})=0$. Therefore, $\psi_{n}$ is $\mathfrak{D}_{0}$-module homomorphism. Since $x\in \mathfrak{D}_{0}, y\in \mathfrak{D}_{n}$, then
\begin{align*}
 \psi_{n}([x, y])&=\varphi_{n}([x, y])-[[x, y], E_{n}]\\&=[\varphi_{n}(x), y]+[x, \varphi_{n}(y)]-[x, [y, E_{n}]]-[y, [E_{n}, x]]\\&=[x, \varphi_{n}(y)]-[x, [y, E_{n}]]\\&=[x, \varphi_{n}(y)-[y, E_{n}]]\\&=[x, \psi_{n}(y)].
\end{align*}
And by virtue of Lemma \ref{RR}, we have ${\rm Hom}_{\mathfrak{D}_{0}}(\mathfrak{D}_{m}, \mathcal{H}_{n})=0$ for $m\neq n, m\neq0$ and $\psi_{n}|_{\mathfrak{D}_{m}}=0$. At last, by virtue of Lemma \ref{PP}, we obtain $\psi_{n}|_{\mathfrak{D}_{0}}=0$.  Furthermore, for $m\in \mathbb{Z}$, $\psi_{n} |_{\mathfrak{D}_{m}}=0$. Therefore $\varphi_{n}\in {\rm Inn}_{\mathbb{C}}(\mathfrak{D}, \mathcal{H})$. Also, by virtue of  Lemma \ref{2.10}, the conclusion is  proved.
\end{proof}

\begin{lemma}\label{2.9}
Let $\mathfrak{D}$ be the mirror Heisenberg-Virasoro algebra and $\mathcal{H}$\ the twisted Heisenberg algebra.
Then $D\in {\rm Der}_{\mathbb{C}}(\mathfrak{D}, \mathcal{H})_{0}$ if and only if $D$ is a linear map from $\mathfrak{D}$ to $\mathcal{H}$ satisfying that
\begin{itemize}
             \item[{\rm{(i)}}] $D(d_{n})=ah_{n+\frac{1}{2}},$
             \item[{\rm{(ii)}}] $D({\bf c})=0,$
             \item[{\rm{(iii)}}] $D(h_{n+\frac{1}{2}})=bh_{n+\frac{1}{2}}+\delta_{n},_{-1}a{\bf l},$
             \item[{\rm{(iv)}}] $D({\bf l})=2b{\bf l},$
           \end{itemize}
           where $n\in \mathbb{Z},\ a,\ b\in \mathbb{C}.$
\end{lemma}
\begin{proof}
Sufficiency: Suppose $D$ is a linear map that satisfies (i)-(iv), then it is easy to verify that $D$ is a derivation from $\mathfrak{D}$ to $\mathcal{H}$.

Necessity:
 (i) Set $D\in{\rm Der}_{\mathbb{C}}(\mathfrak{D}, \mathcal{H})$.
when $n\neq-1$, suppose $$D(d_{n})=f(n)h_{n+\frac{1}{2}}$$ for $m\neq\pm n$, on the one hand
\begin{align*}
D([d_{m}, d_{n}])&=(m-n)D(d_{m+n})\\&=(m-n)f(m+n)h_{m+n+\frac{1}{2}}.
\end{align*}
On the other hand
\begin{align*}
D([d_{m}, d_{n}])&=[D(d_{m}), d_{n}]+[d_{m}, D(d_{n})]\\[5pt]&=f(m)[h_{m+\frac{1}{2}}, d_{n}]+f(n)[d_{m}, h_{n+\frac{1}{2}}]\\&=(m+\frac{1}{2})f(m)h_{m+n+\frac{1}{2}}-(n+\frac{1}{2})f(n)h_{m+n+\frac{1}{2}}.
\end{align*}
Then $$(m-n)f(m+n)=(m+\frac{1}{2})f(m)-(n+\frac{1}{2})f(n).$$
Let $m=1,\ n=0.$ Then $f(0)=f(1)=a$ ($a\in \mathbb{C}$). Next, using induction on $n$, assuming $n-1$, the conclusion holds. Then
\begin{align*}
  (n-2)f(n)&=(n-2)f((n-1)+1)\\&=(n-\frac{1}{2})f(n-1)-\frac{3}{2}f(1)\\&=(n-2)a.
\end{align*}
Hence $f(n)=a$ for $n\in \mathbb{Z}$. When $n = -1$, on the one hand
\begin{align*}
 D([d_{1}, d_{-2}])=3D(d_{-1}).
\end{align*}
On the other hand
\begin{align*}
D([d_{1}, d_{-2}])&=[D(d_{1}), d_{-2}]+[d_{1}, D(d_{-2})]\\&=a[h_{\frac{3}{2}}, d_{-2}]+a[d_{1}, h_{-\frac{3}{2}}]\\&=3ah_{-\frac{1}{2}}.
\end{align*}
Therefore $D(d_{-1})=ah_{-\frac{1}{2}}$.
Thus $D(d_{n})=ah_{n+\frac{1}{2}},\ a\in \mathbb{C}$.

  (ii) On the one hand
\begin{align*}
 D([d_{2}, d_{-2}])&=4D(d_{0})+\frac{1}{2}D({\bf c})\\&=4ah_{\frac{1}{2}}+\frac{1}{2}D({\bf c}).
\end{align*}
  On the other hand
   \begin{align*}
 D([d_{2}, d_{-2}])&=[D(d_{2}), d_{-2}]+[d_{2}, D(d_{-2})]
 \\&=a[h_{\frac{5}{2}}, d_{-2}]+a[d_{2}, h_{-\frac{3}{2}}]
 \\&=4ah_{\frac{1}{2}}.
\end{align*}
Therefore, $D({\bf c})=0.$

 (iii) For $n\neq -1$, let $D(h_{n+\frac{1}{2}})=g(n)h_{n+\frac{1}{2}}$. Since
\begin{align*}
D([d_{n}, h_{m+\frac{1}{2}}])&=(m-\frac{1}{2})D(h_{m+n+\frac{1}{2}})
                           \\&=-(m+\frac{1}{2})g(n+m)h_{m+n+\frac{1}{2}}.
\end{align*}
and
\begin{align*}
D([d_{n}, h_{m+\frac{1}{2}}])&=[D(d_{n}), h_{m+\frac{1}{2}}]+[d_{n}, D(h_{m+\frac{1}{2}})]\\[5pt]&=a[h_{n+\frac{1}{2}}, h_{m+\frac{1}{2}}]+g(m)[d_{n}, h_{m+\frac{1}{2}}]\\&=-(m+\frac{1}{2})g(m)h_{m+n+\frac{1}{2}}.
\end{align*}
Then $g(n+m)=g(m),\ n\neq -1,\ m\neq -n-1, -1$.\\[-5pt]

Set $g(1)=b,\ b\in \mathbb{C}$. Then $g(n)=b$ for $n\in \mathbb{Z}\setminus \{-1\}$.
Also since \[D([d_{-1}, h_{\frac{1}{2}}])=-\frac{1}{2}D(h_{-\frac{1}{2}})\]
and
\begin{align*}
D([d_{-1}, h_{\frac{1}{2}}])&=[D(d_{-1}), h_{\frac{1}{2}}]+[d_{-1}, D(h_{\frac{1}{2}})]\\&=a[h_{-\frac{1}{2}}, h_{\frac{1}{2}}]+b[d_{-1}, h_{\frac{1}{2}}]\\&=-\frac{1}{2}a{\bf l}-\frac{1}{2}bh_{-\frac{1}{2}},
\end{align*}
then $D(h_{-\frac{1}{2}})=a{\bf l}+bh_{-\frac{1}{2}}$.
Thus $D(h_{n+\frac{1}{2}})=bh_{n+\frac{1}{2}}+\delta_{n,-1}a{\bf l}$.

  (iv) Since
  \[
  D([h_{-\frac{1}{2}}, h_{\frac{1}{2}}])=-\frac{1}{2}D({\bf l})
  \]
  and
  \begin{align*}
D([h_{-\frac{1}{2}}, h_{\frac{1}{2}}])&=[D(h_{-\frac{1}{2}}), h_{\frac{1}{2}}]+[h_{-\frac{1}{2}}, D(h_{\frac{1}{2}}]\\&=b[h_{-\frac{1}{2}}, h_{\frac{1}{2}}]+b[h_{-\frac{1}{2}}, h_{\frac{1}{2}}]\\&=-b{\bf l},
\end{align*}
then $D({\bf l})=2b{\bf l}$.
\end{proof}

According to Lemma \ref{2.9}, $D_{i}$ $(i=1,2)$ is given by:
\begin{align*}
&D_1\colon D_{1}|_{\mathcal{V}}=0,\ D_{1}(h_{n+\frac12})=h_{n+\frac12},\ D_{{\bf 1}}({\bf l})=2{\bf l}\ (n\in\mathbb{Z}),\\ 
&D_2\colon D_{2}(d_{n})=h_{n+\frac{1}{2}},\ D_{2}(h_{n+\frac{1}{2}})=\delta_{n,-1}l,\ D_{2}({\bf c})=D_{2}({\bf l})=0\,(n\in\mathbb{Z}).
\end{align*}

\begin{lemma}\label{PO}
Let $\mathfrak{D}$ be the mirror Heisenberg-Virasoro algebra and $\mathcal{H}$\ the twisted Heisenberg algebra. Then
$${\rm H}^{1}(\mathfrak{D}, \mathcal{H})=\mathbb{C}D_{1}\oplus \mathbb{C}D_{2}.$$
\begin{proof}
Let $$k_{1}D_{1}+k_{2}D_{2}, k_{i}\in \mathbb{C}, i=1,2.$$ Then for $h_{-\frac12}$
\[
0=(k_{1}D_{1}+k_{2}D_{2})({h_{-\frac{1}{2}}})
=k_{1}h_{-\frac{1}{2}}+k_{2}l.
\]
Therefore, $k_{1}=k_{2}=0$. That is, $\{D_{i}\mid i=1,2\}$ are linearly independent. And for any $D\in {\rm Der}_{\mathbb{C}}(\mathfrak{D}, \mathcal{H})_{0}$ (Because of the linearity of the derivation, it suffices to verify the action on the basis), we have
\begin{align*}
  &D(d_n)=ah_{n+\frac{1}{2}}=aD_2(d_{n}),\\
  &D({\bf c})=D_1({\bf c})=0,\\
  &D(h_{n+\frac{1}{2}})=bh_{n+\frac{1}{2}}+\delta_{n},_{-1}a{\bf l}
  =bD_1(h_{n+\frac{1}{2}})+aD_2(h_{n+\frac{1}{2}}),\\
  &D({\bf l})=2b{\bf l}=bD_1({\bf l}).
\end{align*}
By virtue of Lemma \ref{QQ}, we obtain
$${\rm Der}_{\mathbb{C}}(\mathfrak{D}, \mathcal{H})=\mathbb{C}D_{1}\oplus \mathbb{C}D_{2}\oplus{\rm Inn}_{\mathbb{C}}(\mathfrak{D}, \mathcal{H})$$
and  $${\rm H}^{1}(\mathfrak{D}, \mathcal{H})={\rm Der}_{\mathbb{C}}(\mathfrak{D}, \mathcal{H})/ {\rm Inn}_{\mathbb{C}}(\mathfrak{D}, \mathcal{H}).$$
\end{proof}
\end{lemma}

\begin{lemma}\label{ZA}
 Let\ $\mathcal{V}$ be the Virasoro algebra and $\mathcal{H}$\,the twisted Heisenberg algebra. Then $${\rm Hom}_{{\mathcal{U}}(\mathcal{V})}(\mathcal{H} /[\mathcal{H}, \mathcal{H}], \mathcal{V}]=0.$$
\end{lemma}
\begin{proof}
By the definition of $\mathcal{H}$, we have $[\mathcal{H}, \mathcal{H}]=\mathbb{C}{\bf l}$. Put $f\in {\rm Hom}_{{\mathcal{U}}(\mathcal{V})}(\mathcal{H} / [\mathcal{H}, \mathcal{H}], \mathcal{V}]$, for $m\in \mathbb{Z}$.
Suppose \[f(h_{m+\frac{1}{2}})=\sum a_{i}(m)d_{m_{i}}+k{\bf c},\, k\in \mathbb{C},\]  on the one hand
\begin{align*}
  [d_{0},f(h_{m+\frac{1}{2}})]&=[d_{0},\sum a_{i}(m)d_{m_{i}}+k{\bf c}]
  \\&=\sum a_{i}(m)[d_{0}, d_{m_{i}}]
  \\&=-\sum a_{i}(m)m_{i}d_{m_{i}},
\end{align*}
on the other hand
\begin{align*}
[d_{0},f(h_{m+\frac{1}{2}})]&=f([d_{0}, h_{m+\frac{1}{2}}])\\&=-(m+\frac{1}{2})f(h_{m+\frac{1}{2}})\\&=-(m+\frac{1}{2})(\sum a_{i}(m)d_{m_{i}}+k{\bf c})\\&=-(m+\frac{1}{2})\sum a_{i}(m)d_{m_{i}}-(m+\frac{1}{2})k{\bf c},
\end{align*}
then \begin{align*} a_{i}(m)m_{i}&=(m+\frac{1}{2})a_{i}(m),\\ k&=0.\end{align*}
Since $m_{i}\in\mathbb{Z},\ m\in\mathbb{Z}$, $m_{i}\neq m+\frac{1}{2}$, then $a_{i}(m)=0$.
Furthermore, we obtain ${\rm Hom}_{{\mathcal{U}}(\mathcal{V})}(\mathcal{H} / [\mathcal{H}, \mathcal{H}], \mathcal{V}]=0$.
\end{proof}

\begin{theorem}\label{VB}
Let $\mathfrak{D}$ be the mirror Heisenberg-Virasoro algebra and $\mathcal{H}$\ the twisted Heisenberg algebra. Then
$${\rm Der}_{\mathbb{C}}\mathfrak{D}=\mathbb{C}D_{1}\oplus \mathbb{C}D_{2}\oplus {\rm ad}\mathfrak{D}.$$
\begin{proof}
Since $${\rm H}^{1}(\mathcal{H}, \mathfrak{D}/\mathcal{H})^{\mathfrak{D}/\mathcal{H}}
 \cong {\rm Hom}_{\mathbb{C}}(\mathcal{H}/[\mathcal{H}, \mathcal{H}], \mathfrak{D}/\mathcal{H})^{\mathfrak{D}/\mathcal{H}}
 ={\rm Hom}_{\mathfrak{D}/\mathcal{H}}(\mathcal{H}/[\mathcal{H}, \mathcal{H}], \mathfrak{D}/\mathcal{H}),$$
by virtue of Lemma \ref{ZA}, then ${\rm H}^{1}(\mathcal{H}, \mathfrak{D}/\mathcal{H})^{\mathfrak{D}/\mathcal{H}}=0$. According to the exact sequence (\ref{2})
 \begin{equation*}
 0 \rightarrow {\rm H}^{1}(\mathfrak{D}/ \mathcal{H}, \mathfrak{D}/ \mathcal{H})\rightarrow {\rm H}^{1}(\mathfrak{D}, \mathfrak{D}/\mathcal{H})\rightarrow {\rm H}^{1}(\mathcal{H}, \mathfrak{D}/\mathcal{H})^{\mathfrak{D}/\mathcal{H}}
 \end{equation*}
 and ${\rm H}^{1}(\mathfrak{D}/\mathcal{H}, \mathfrak{D}/\mathcal{H})=0$, we obtain ${\rm H}^{1}(\mathfrak{D}, \mathfrak{D}/\mathcal{H})=0$. Then by virtue of Lemma \ref{PO} and the exact sequence (\ref{9})
 \begin{equation*}
{\rm H}^{1}(\mathfrak{D}, \mathcal{H})\rightarrow {\rm H}^{1}(\mathfrak{D}, \mathfrak{D})\rightarrow {\rm H}^{1}(\mathfrak{D}, \mathfrak{D}/\mathcal{H}) \label{1},
\end{equation*}
we have  $${\rm H}^{1}(\mathfrak{D}, \mathcal{H})={\rm H}^{1}(\mathfrak{D}, \mathfrak{D}).$$
Therefore
$${\rm H}^{1}(\mathfrak{D}, \mathfrak{D})=\mathbb{C}D_{1}\oplus \mathbb{C}D_{2}.$$
Since ${\rm H}^{1}(\mathfrak{D}, \mathfrak{D})={\rm Der}\mathfrak{D}/{\rm ad}\mathfrak{D}$, then the conclusion holds.
\end{proof}
\end{theorem}

\section{2-local derivations on the Mirror Heisenberg-Virasoro algebra}

Let $L$ be an algebra. $\Delta$ is called a 2-local derivation of $L$, if for every $ x, y\in L$, there exists a derivation $\Delta_{x, y}$ such that $$\Delta_{x, y}(x)=\Delta(x), \Delta_{x, y}(y)=\Delta(y).$$

Let $\Delta$ be a 2-local derivation on the mirror Heisenberg-Virasoro algebra $\mathfrak{D}$. Then for any $ x, y\in \mathfrak{D} $, there exists $\Delta_{x, y}$ satisfying $$\Delta_{x, y}(x)=\Delta(x), \Delta_{x, y}(y)=\Delta(y).$$ By Theorem \ref{VB}, it can be written as
\begin{align*}
\Delta_{x, y}=&{\rm{ad}}(\sum_{i\in \mathbb{Z}}a_{i}(x, y)d_{i}+ b_{i}(x, y)h_{i+\frac{1}{2}}+ l_{1}(x, y){\bf c}+ l_{2}(x, y){\bf l} )
 \\&+\alpha(x, y)D_{1}+\beta(x, y)D_{2}\tag{3.1}\label{LO}.
\end{align*}
where $a_{i},\ b_{i}(i\in \mathbb{Z}),\ l_{1},\ l_{2},\ \alpha,\ \beta$ are complex-valued functions on $\mathfrak{D} \times \mathfrak{D}$ and  $D_{k}$ is given in Theorem \ref{VB}, where $k=1,2.$

\begin{lemma}\label{RF}
 Let $\Delta$ be a 2-local derivation on the mirror Heisenberg-Virasoro algebra $\mathfrak{D}$. Then for any fixed $x\in\mathfrak{D}$,
\begin{itemize}
  \item[{\rm{(i)}}]
   If $ \Delta(d_{i})=0$ for  $ i\in \mathbb{Z} $, then
\begin{align*}
\Delta_{d_{i}, x}=&{\rm ad}\big(a_{i}(d_{i}, x)d_{i}+ b_{0}(d_{i}, x)h_{\frac{1}{2}}+l_{1}(d_{i}, x){\bf c}+l_{2}(d_{i}, x){\bf l} \big)
                \\&+\alpha(d_{i}, x)D_{1}-\frac{1}{2}b_{0}(d_{i}, x)D_{2}.
\end{align*}
  \item[{\rm{(ii)}}]
 If $\Delta(h_{\frac{1}{2}})=0$, then
\begin{align*}
\Delta_{h_{\frac{1}{2}}, x}=&{\rm ad}\big(a_{0}(h_{\frac{1}{2}}, x)d_{0}+ \sum_{i\neq-1} b_{i}(h_{\frac{1}{2}}, x)h_{i+\frac{1}{2}}+ l_{1}(h_{\frac{1}{2}}, x){\bf c}+ l_{2}(h_{\frac{1}{2}}, x){\bf l}\big)\\&+\frac{1}{2}a_{0}(h_{\frac{1}{2}}, x)D_{1}+\beta(h_{\frac{1}{2}}, x)D_{2}.
\end{align*}
  \item[{\rm{(iii)}}]
 If $\Delta(h_{-\frac{1}{2}})=0$, then
\begin{align*}
\Delta_{h_{-\frac{1}{2}}, x}=&{\rm ad}\big(a_{0}(h_{-\frac{1}{2}}, x)d_{0}+ \sum_{i\in\mathbb{Z}} b_{i}(h_{-\frac{1}{2}}, x)h_{i+\frac{1}{2}}+ l_{1}(h_{\frac{1}{2}}, x){\bf c}+ l_{2}(h_{-\frac{1}{2}}, x){\bf l}\big)\\&-\frac{1}{2}a_{0}(h_{-\frac{1}{2}}, x)D_{1}-\frac{1}{2}b_0(h_{-\frac{1}{2}}, x)D_{2}.
\end{align*}
  \item[{\rm{(iv)}}]
  If $\Delta(h_{i+\frac{1}{2}})=0$ for $i\neq0,-1$, then
\begin{align*}
\Delta_{h_{i+\frac{1}{2}}, x}=&{\rm ad}\big(a_{0}(h_{i+\frac{1}{2}}, x)d_{0}+ \sum_{j\neq-1-i}b_{j}(h_{i+\frac{1}{2}}, x)h_{j+\frac{1}{2}}+ l_{1}(h_{i+\frac{1}{2}}, x){\bf c}+ l_{2}(h_{i+\frac{1}{2}}, x){\bf l} \big)\\&+(i+\frac{1}{2})a_{0}(h_{i+\frac{1}{2}}, x)D_{1}+\beta(h_{i+\frac{1}{2}}, x)D_{2}.
\end{align*}
\end{itemize}
\end{lemma}
\begin{proof}
\begin{itemize}
  \item[(i)]
 According to $\Delta(d_{i})=0$ and the formula (\ref{LO}), we have
  \begin{alignat*}{2}
\Delta(d_{i})&=& &\Delta_{d_{i}, x}(d_{i})
\\[5pt]&=&& \big[ \sum_{j\in \mathbb{Z}}a_{j}(d_{i}, x)d_{j}+ b_{j}(d_{i}, x)h_{j+\frac{1}{2}}+ l_{1}(d_{i}, x){\bf c}+ l_{2}(d_{i}, x){\bf l}, d_{i}\big]
\\ &&&+\alpha(d_{i}, x)D_{1}(d_{i})+\beta(d_{i}, x)D_{2}(d_{i})
\\ &=&&\sum_{j\in\mathbb{Z}}(j-i)a_{j}(d_{i}, x)d_{i+j}+(j+\frac{1}{2})b_{j}(d_{i}, x)h_{i+j+\frac{1}{2}}
 \\&&&-\frac{i^{3}-i}{12}a_{-i}(d_{i}, x)c+\beta(d_{i}, x)h_{i+\frac{1}{2}}
 \\&=&&0.
\end{alignat*}
Then we have
\begin{align*}
  &a_{j}(d_{i}, x)=0 \,\, \mbox{for}\,j\neq i,\\[5pt]
  &b_{j}(d_{i}, x)=0 \,\, \mbox{for}\,  j\neq0,\\
 &\beta(d_{i}, x)=-\frac{1}{2}b_{0}(d_{i}, x).
\end{align*}
\item[(ii)]
According to $\Delta(h_{\frac{1}{2}})=0$ and the formula (\ref{LO}), we have
\begin{alignat*}{2}
\Delta(h_{\frac{1}{2}})&=& &\Delta_{h_{\frac{1}{2}}, x}(h_{\frac{1}{2}})
\\&=& &\big[ \sum_{j\in \mathbb{Z}}a_{j}(h_{\frac{1}{2}}, x)d_{j}+ b_{j}(h_{\frac{1}{2},}, x)h_{j+\frac{1}{2}}+ l_{1}(h_{\frac{1}{2}}, x){\bf c}+ l_{2}h_{\frac{1}{2}, x}{\bf l}, h_{\frac{1}{2}}\big]
\\& & &+\alpha(h_{\frac{1}{2}}, x)D_{1}(h_{\frac{1}{2}})+\beta(h_{\frac{1}{2}}, x)D_{2}(h_{\frac{1}{2}})
\\&=& &\sum_{j\in\mathbb{Z}}(-\frac{1}{2})a_{j}(h_{\frac{1}{2}}, x)h_{j+\frac{1}{2}}-\frac{1}{2}b_{-1}(h_{\frac{1}{2}}, x){\bf l}+\alpha(h_{\frac{1}{2}}, x)h_{\frac{1}{2}}
\\&=& &0.
\end{alignat*}
Then we have
\begin{align*}
  &a_{j}(h_{\frac{1}{2}}, x)=0\,\, \mbox{for}\,  j\neq 0,\\[5pt]
   &b_{-1}(h_{\frac{1}{2}}, x)=0,\\
&\alpha(h_{\frac{1}{2}}, x)=\frac{1}{2}a_{0}(h_{\frac{1}{2}}, x).
\end{align*}
  \item[(iii)]
 By $\Delta(h_{-\frac{1}{2}})=0$  and the formula (\ref{LO}), we have
  \begin{alignat*}{2}
\Delta(h_{-\frac{1}{2}})&=&&\Delta_{h_{-\frac{1}{2}}, x}(h_{-\frac{1}{2}})
\\&=&&\big[\sum_{j\in \mathbb{Z}}a_{j}(h_{-\frac{1}{2}}, x)d_{j}+ b_{j}(h_{-\frac{1}{2}}, x)h_{j+\frac{1}{2}}+ l_{1}(h_{-\frac{1}{2}}, x){\bf c}+ l_{2}(h_{-\frac{1}{2}}, x){\bf l}, h_{-\frac{1}{2}}\big]
\\&&&+\alpha(h_{-\frac{1}{2}}, x)D_{1}(h_{-\frac{1}{2}})+\beta(h_{-\frac{1}{2}}, x)D_{2}(h_{-\frac{1}{2}})
\\&=&&\sum_{j\in\mathbb{Z}}\frac{1}{2}a_{j}(h_{-\frac{1}{2}}, x)h_{j-\frac{1}{2}}+\frac{1}{2}b_{0}(h_{-\frac{1}{2}}, x){\bf l}+\alpha(h_{-\frac{1}{2}}, x)h_{-\frac{1}{2}}+\beta (h_{-\frac{1}{2}}, x){\bf l}
\\&=&&0.
\end{alignat*}
Then we have
\begin{align*}
&a_{j}(h_{-\frac{1}{2}}, x)=0\,\,\mbox{for}\,\,  j\neq 0,\\[5pt]
&\alpha(h_{-\frac{1}{2}}, x)=-\frac{1}{2}a_{0}(h_{-\frac{1}{2}}, x),\\
&\beta(h_{-\frac{1}{2}}, x)=-\frac{1}{2}b_{0}(h_{-\frac{1}{2}}, x).
\end{align*}
 \item[(iv)]
 By $\Delta(h_{i+\frac{1}{2}})=0$ and the formula (\ref{LO}), we have
 \begin{alignat*}{2}
\Delta(h_{i+\frac{1}{2}})&=&&\Delta_{h_{i+\frac{1}{2}}, x}(h_{i+\frac{1}{2}})
\\[5pt]&=&&\big[\sum_{j\in \mathbb{Z}}a_{j}(h_{i+\frac{1}{2}}, x)d_{j}+ b_{j}(h_{i+\frac{1}{2}}, x)h_{j+\frac{1}{2}}+ l_{1}(h_{i+\frac{1}{2}}, x){\bf c}+ l_{2}(h_{i+\frac{1}{2}}, x){\bf l}, h_{i+\frac{1}{2}}\big]
\\&&&+\alpha(h_{i+\frac{1}{2}}, x)D_{1}(h_{i+\frac{1}{2}})+\beta(h_{i+\frac{1}{2}}, x)D_{2}(h_{i+\frac{1}{2}})
\\&=&&\sum_{j\in\mathbb{Z}}-(i+\frac{1}{2})a_{j}(h_{i+\frac{1}{2}}, x)h_{j+i+\frac{1}{2}}+(-i-\frac{1}{2})b_{-1-i}(h_{i+\frac{1}{2}}, x){\bf l}+\alpha(h_{i+\frac{1}{2}}, x)h_{i+\frac{1}{2}}\\&=&&0.
\end{alignat*}
Then we have
\begin{align*}
 &b_{j}(h_{i+\frac{1}{2}}, x)=0\,\,\mbox{for}\,\, j=-1-i,\\[5pt]
  &a_{j}(h_{i+\frac{1}{2}}, x)=0\,\,\mbox{for}\,\, j\neq0,\\
&\alpha(h_{i+\frac{1}{2}}, x)=(i+\frac{1}{2})a_{0}(h_{i+\frac{1}{2}}, x).
 \end{align*}
\end{itemize}
\end{proof}

\begin{lemma}\label{KM}
Let $\Delta$ be a 2-local derivation on the mirror Heisenberg-Virasoro algebra $\mathfrak{D}$ such that $\Delta(d_{0})=\Delta(d_{1})=0$. Then $\Delta(d_{i})=0$ for  $i\in\mathbb{Z}$.
\end{lemma}
\begin{proof}
Put $k=0,1$. By Lemma \ref{RF}, suppose
\begin{align*}
\Delta_{d_{k}, x}=&{\rm ad} a_{k}(d_{k}, x)d_{k}+b_{0}(d_{k}, x)h_{\frac{1}{2}}+l_{1}(d_{k}, x){\bf c}+l_{2}(d_{k}, x){\bf l}\\&+\alpha(d_{k}, x)D_{1}-\frac{1}{2}b_{0}(d_{k}, x)D_{2},
\end{align*}
then
\begin{alignat*}{2}
\Delta(d_{i})&=&&\Delta_{d_{0}, d_{i}}(d_{i})\\[5pt]&=&&[a_{0}(d_{0}, d_{i})d_{0}+b_{0}(d_{0}, d_{i})h_{\frac{1}{2}}+l_{1}(d_{0}, d_{i}){\bf c}+l_{2}(d_{0}, d_{i}){\bf l}, d_{i}]
\\&&&+\alpha(d_{0}, d_{i})D_{1}(d_{i})-\frac{1}{2}b_{0}(d_{0}, d_{i})D_{2}(d_{i})
\\&=&&-ia_{0}(d_{0}, d_{i})d_{i}+\frac{1}{2}b_{0}(d_{0}, d_{i})h_{i+\frac{1}{2}}
-\frac{1}{2}b_{0}(d_{0}, d_{i})h_{i+\frac{1}{2}}\\&=&&-ia_{0}(d_{0}, d_{i})d_{i}.
\end{alignat*}
Another,
\begin{alignat*}{2}
\Delta(d_{i})&=&&\Delta_{d_{1}, d_{i}}(d_{i})\\[5pt]&=&&[a_{1}(d_{1}, d_{i})d_{1}+b_{0}(d_{1}, d_{i})h_{\frac{1}{2}}+l_{1}(d_{1}, d_{i}){\bf c}+l_{2}(d_{1}, d_{i}){\bf l}, d_{i}]\\&&&+\alpha(d_{1}, d_{i})D_{1}(d_{i})-\frac{1}{2}b_{0}(d_{1}, d_{i})D_{2}(d_{i})\\&=&&(1-i)a_{1}(d_{1}, d_{i})d_{i+1}+\frac{1}{2}b_{0}(d_{1}, d_{i})h_{i+\frac{1}{2}}-\frac{1}{2}b_{0}(d_{1}, d_{i})h_{i+\frac{1}{2}}\\&=&&(1-i)a_{1}(d_{1}, d_{i})d_{i+1}.
\end{alignat*}
Then we have $$a_{0}(d_{0}, d_{i})=a_{1}(d_{1}, d_{i})=0.$$
Therefore $\Delta(d_{i})=0$ for  $i\in\mathbb{Z}$.
\end{proof}

\begin{lemma}\label{3.4}
Let $\Delta$ be a 2-local derivation on the mirror Heisenberg-Virasoro algebra $\mathfrak{D}$ such that $\Delta(d_{i})=0$. Then
$$\Delta(x)=\lambda_{x}(\sum_{t\in\mathbb{Z}}\beta_{t}h_{t+\frac{1}{2}}+
2k_{2}{\bf l})$$ for any $$x=\sum_{t\in\mathbb{Z}}(\alpha_{t}d_{t}+\beta_{t}h_{t+\frac{1}{2}})+k_{1}{\bf c}+k_{2}{\bf l}\in\mathfrak{D},$$
where $\lambda_x$ is complex number depending on $x$.
\end{lemma}
\begin{proof}
For $x=\sum_{t\in\mathbb{Z}}(\alpha_{t}d_{t}+\beta_{t}h_{t+\frac{1}{2}})+k_{1}{\bf c}+k_{2}{\bf l}\in\mathfrak{D}$, since $\Delta(d_{i})=0$, by Lemma \\[5pt]\ref{RF}, then we have
\begin{alignat*}{2}
\Delta(x)&=&&\Delta_{d_{0}, x}(x)\\[5pt]&=&&\big({\rm ad}(a_{0}(d_{0}, x)d_{0}+b_{0}(d_{0}, x)h_{\frac{1}{2}}+l_{1}(d_{0}, x){\bf c}+l_{2}(d_{0}, x){\bf l})\\&&&+\alpha(d_{0}, x)D_{1}-\frac{1}{2}b_{0}(d_{0}, x)D_{2}\big)(x)\\&=&&a_{0}(d_{0}, x)(\sum_{t\in\mathbb{Z}}-t\alpha_{t}d_{t}-(t+\frac{1}{2})\beta_{t}h_{t+\frac{1}{2}})+b_{0}(d_{0}, x)(\sum_{t\in\mathbb{Z}}\frac{1}{2}\alpha_{t}h_{t+\frac{1}{2}}+\frac{1}{2}\beta_{-1}{\bf l})\\&&&+\alpha(d_{0}, x)(\sum_{t\in\mathbb{Z}}\beta_{t}h_{t+\frac{1}{2}}+2k_{2}{\bf l})-\frac{1}{2}b_{0}(d_{0}, x)(\sum_{t\in\mathbb{Z}}\alpha_{t}h_{t+\frac{1}{2}}+\beta_{-1}{\bf l})\\&=&&a_{0}(d_{0}, x)(\sum_{t\in\mathbb{Z}}-t\alpha_{t}d_{t}-(t+\frac{1}{2})\beta_{t}h_{t+\frac{1}{2}})+\alpha(d_{0}, x)(\sum_{t\in\mathbb{Z}}\beta_{t}h_{t+\frac{1}{2}}+2k_{2}{\bf l}).
\end{alignat*}
Another,
\begin{alignat*}{2}
\Delta(x)&=&&\Delta_{d_{i}, x}(x)\\[5pt]&=&&\big({\rm{ad}}(a_{i}(d_{i}, x)d_{i}+b_{0}(d_{i}, x)h_{\frac{1}{2}}+l_{1}(d_{i}, x){\bf c}+l_{2}(d_{i}, x){\bf l})\\&&&+\alpha(d_{i}, x)D_{1}-\frac{1}{2}b_{0}(d_{i}, x)D_{2}\big)(x)\\&=&&a_{i}(d_{i}, x)(\sum_{t\in\mathbb{Z}}(i-t)\alpha_{t}d_{i+t}-(t+\frac{1}{2})\beta_{t}h_{i+t+\frac{1}{2}})+b_{0}(d_{i}, x)(\sum_{t\in\mathbb{Z}}\frac{1}{2}\alpha_{t}h_{t+\frac{1}{2}}+\frac{1}{2}\beta_{-1}{\bf l})
\\&&&+\alpha(d_{i}, x)(\sum_{t\in\mathbb{Z}}\beta_{t}h_{t+\frac{1}{2}}+2k_{2}{\bf l})-\frac{1}{2}b_{0}(d_{i}, x)(\sum_{t\in\mathbb{Z}}\alpha_{t}h_{t+\frac{1}{2}}+\beta_{-1}{\bf l})-\frac{i^{3}-i}{12}a_{-i}(d_{i}, x){\bf c}\\&=&&a_{i}(d_{i}, x)\big(\sum_{t\in\mathbb{Z}}(i-t)\alpha_{t}d_{i+t}-(t+\frac{1}{2})\beta_{t}h_{i+t+\frac{1}{2}}\big)+\alpha(d_{i}, x)\big(\sum_{t\in\mathbb{Z}}\beta_{t}h_{t+\frac{1}{2}}+2k_{2}{\bf l}\big)\\&&&-\frac{i^{3}-i}{12}a_{-i}(d_{i}, x){\bf c}.
\end{alignat*}
Comparing the two equations above,
when neither $\alpha _{t}$ nor $\beta_{t}$ is equal to 0, by taking different $i\in\mathbb{Z}$, we have
\begin{align*}
a_{0}(d_{0}, x)&=a_{i}(d_{i}, x)=0,\\
 \alpha(d_{0}, x)&=\alpha(d_{i}, x).
\end{align*}
 If $\alpha _{t}=0, \beta_{t}\neq0$, by taking different $i\in\mathbb{Z}$, then we have
\begin{align*}
a_{0}(d_{0}, x)&=a_{i}(d_{i}, x)=0,\\
 \alpha(d_{0}, x)&=\alpha(d_{i}, x).
\end{align*}
If $\alpha _{t}\neq0, \beta_{t}=0$, then we have
\begin{align*}
a_{0}(d_{0}, x)&=a_{i}(d_{i}, x)=0,\\
 \alpha(d_{0}, x)&=\alpha(d_{i}, x).
\end{align*}
If $\alpha _{t}=\beta_{t}=0$, then we have \[ \alpha(d_{0}, x)=\alpha(d_{i}, x).\]
Thus we have
\[
\Delta(x)=\lambda_{x}(\sum_{t\in\mathbb{Z}}\beta_{t}h_{t+\frac{1}{2}}
+2k_{2}{\bf l}),\]
where $\lambda_{x}=\alpha(d_{i}, x)$  are complex-valued numbers depending on $x$.
\end{proof}

\begin{lemma}\label{3.5}
Let $\Delta$ be a 2-local derivation on the mirror Heisenberg-Virasoro algebra $\mathfrak{D}$ such that $\Delta(d_{0})=\Delta(d_{1})=0$ and there exists $\Delta( h_{t-\frac{1}{2}})=0$ for fixed $t\in\mathbb{Z}\!\setminus\!\{0, 1\}$. Then $\Delta(x)=0$ for $x\in \mathfrak{D}$.
\end{lemma}
\begin{proof}
We proceed by two steps:

${\rm{(i)}}$  We  prove $x=d_{2t}+h_{t+\frac{1}{2}},\ t\in\mathbb{Z}\backslash\{0,1\}$, then
\begin{align*}
\Delta_{d_{2t}+h_{t+\frac{1}{2}}, y}=&{\rm ad}(a_{2t}(d_{2t}+h_{t+\frac{1}{2}}, y)(d_{2t}+h_{t+\frac{1}{2}})+b_{0}(d_{2t}+h_{t+\frac{1}{2}}, y)h_{\frac{1}{2}}\\&+b_{-t}(d_{2t}+h_{t+\frac{1}{2}}, y)h_{-t+\frac{1}{2}}+l_{1}(d_{2t}+h_{t+\frac{1}{2}},y){\bf c}+l_{2}(d_{2t}+h_{t+\frac{1}{2}}, y){\bf l})\\&+(t-\frac{1}{2})b_{-t}(d_{2t}+h_{t+\frac{1}{2}}, y)D_{1}-\frac{1}{2}b_{0}(d_{2t}+h_{t+\frac{1}{2}}, y)D_{2}.
\end{align*}

On the one hand, according to Lemma \ref{KM},
 \[
\Delta(d_{2t}+h_{t+\frac{1}{2}})=\lambda_{d_{2t}+h_{t+\frac{1}{2}}}h_{t+\frac{1}{2}}.
\]
On the other hand, since $\Delta(h_{t-\frac{1}{2}})=0$, then
\begin{alignat*}{2}
\Delta(d_{2t}+h_{t+\frac{1}{2}})&=&&\Delta_{h_{t-\frac{1}{2}}, d_{2t}+h_{t+\frac{1}{2}}}(d_{2t}+h_{t+\frac{1}{2}})
\\[5pt]&=&&\big({\rm ad}(a_{0}(h_{t-\frac{1}{2}}, d_{2t}+h_{t+\frac{1}{2}})d_{0}
+\sum_{j\neq-t}b_{j}(h_{t-\frac{1}{2}}, d_{2t}+h_{t+\frac{1}{2}})h_{j+\frac{1}{2}}\\&&&+l_{1}(h_{t-\frac{1}{2}}, d_{2t}+h_{t+\frac{1}{2}}){\bf c}+l_{2}(h_{t-\frac{1}{2}}, d_{2t}+h_{t+\frac{1}{2}}){\bf l})\\&&&+(t-\frac{1}{2})a_{0}(h_{t-\frac{1}{2}}, d_{2t}+h_{t+\frac{1}{2}})D_{1}+\beta(h_{t-\frac{1}{2}}, d_{2t}+h_{t+\frac{1}{2}})D_{2}\big)(d_{2t}+h_{t+\frac{1}{2}})
\\&=&&a_{0}(h_{t-\frac{1}{2}}, d_{2t}+h_{t+\frac{1}{2}})(-2td_{2t}-(t+\frac{1}{2})h_{t+\frac{1}{2}})
\\&&&+\sum_{j\neq-t}b_{j}(h_{t-\frac{1}{2}}, d_{2t}+h_{t+\frac{1}{2}})(j+\frac{1}{2})h_{2t+j+\frac{1}{2}}\\&&&+(-t-\frac{1}{2})b_{-1-t}(h_{t-\frac{1}{2}}, d_{2t}+h_{t+\frac{1}{2}}){\bf l}
\\&&&+(t-\frac{1}{2})a_{0}(h_{t-\frac{1}{2}}, d_{2t}+h_{t+\frac{1}{2}})h_{t+\frac{1}{2}}+\beta(h_{t-\frac{1}{2}}, d_{2t}+h_{t+\frac{1}{2}})(h_{2t+\frac{1}{2}}+\delta_{t, -1}{\bf l}).
\end{alignat*}
Comparing the two equations above, we have
\begin{align*}
&a_{0}(h_{t-\frac{1}{2}}, d_{2t}+h_{t+\frac{1}{2}})=0,\\[5pt]
&b_{-1-t}(h_{t-\frac{1}{2}}, d_{2t}+h_{t+\frac{1}{2}})=0, t\neq-1,\\
&\beta(h_{t-\frac{1}{2}}, d_{2t}+h_{t+\frac{1}{2}})=-\frac{1}{2}b_{0}(h_{t-\frac{1}{2}}, d_{2t}+h_{t+\frac{1}{2}}).
\end{align*}
Therefore, $\lambda_{d_{2t}+h_{t+\frac{1}{2}}}=0$. Furthermore, we have $\Delta(d_{2t}+h_{t+\frac{1}{2}})=0.$

By Lemma \ref{LO}, set
\begin{alignat*}{2}
\Delta_{{d_{2t}+h_{t+\frac{1}{2}}}, y}&=&&{\rm{ad}}(\sum_{i\in \mathbb{Z}}a_{i}(d_{2t}+h_{t+\frac{1}{2}}, y)d_{i}+ b_{i}(d_{2t}+h_{t+\frac{1}{2}}, y)h_{i+\frac{1}{2}}
\\&&&+ l_{1}(d_{2t}+h_{t+\frac{1}{2}}, y){\bf c}+ l_{2}(d_{2t}+h_{t+\frac{1}{2}}, y){\bf l} )\\[5pt]&&&+\alpha(d_{2t}+h_{t+\frac{1}{2}}, y)D_{1}
+\beta(d_{2t}+h_{t+\frac{1}{2}}, y)D_{2}\\[5pt]
\Delta(d_{2t}+h_{t+\frac{1}{2}})&=&&\Delta_{d_{2t}+h_{t+\frac{1}{2}}, y}(d_{2t}+h_{t+\frac{1}{2}})\\[5pt]&=&&\big({\rm{ad}}(\sum_{i\in \mathbb{Z}}a_{i}(d_{2t}+h_{t+\frac{1}{2}}, y)d_{i}+ b_{i}(d_{2t}+h_{t+\frac{1}{2}}, y)h_{i+\frac{1}{2}}
\\[5pt]&&&+ l_{1}(d_{2t}+h_{t+\frac{1}{2}}, y){\bf c}+ l_{2}(d_{2t}+h_{t+\frac{1}{2}}, y){\bf l} )\\[5pt]&&&+\alpha(d_{2t}+h_{t+\frac{1}{2}}, y)D_{1}
+\beta(d_{2t}+h_{t+\frac{1}{2}}, y)D_{2}\big)(d_{2t}+h_{t+\frac{1}{2}})
\\[5pt]&=&&\sum_{i\in \mathbb{Z}}a_{i}(d_{2t}+h_{t+\frac{1}{2}}, y)\big((i-2t)d_{2t+i}-(t+\frac{1}{2})h_{i+t+\frac{1}{2}}\big)
\\&&&-\frac{4t^{3}-t}{6}a_{-2t}(d_{2t}+h_{t+\frac{1}{2}}, y){\bf c}+\sum_{i\in \mathbb{Z}}b_{i}(d_{2t}+h_{t+\frac{1}{2}}, y)(i+\frac{1}{2})h_{2t+i+\frac{1}{2}}
\\&&&+(-t-\frac{1}{2})b_{-1-t}(d_{2t}+h_{t+\frac{1}{2}}, y){\bf l}+\alpha(h_{d_{2t}+h_{t+\frac{1}{2}}}, y)h_{t+\frac{1}{2}}\\[5pt]&&&+\beta(d_{2t}+h_{t+\frac{1}{2}}, y)(h_{2t+\frac{1}{2}}
+\delta_{t,-1}{\bf l})\\&=&&0.
\end{alignat*}
So we have
\begin{align*}
&a_{i}(d_{2t}+h_{t+\frac{1}{2}})=0\ (i\neq2t),\\
&\alpha(h_{d_{2t}+h_{t+\frac{1}{2}}}, y)=(t-\frac{1}{2} )b_{-t}(d_{2t}+h_{t+\frac{1}{2}}),\\ &\beta(d_{2t}+h_{t+\frac{1}{2}}, y)=-\frac{1}{2}b_{0}(d_{2t}+h_{t+\frac{1}{2}}, y), \\ &a_{2t}(d_{2t}+h_{t+\frac{1}{2}}, y)=b_{t}(d_{2t}+h_{t+\frac{1}{2}}, y), \\[5pt]
&b_{i}(d_{2t}+h_{t+\frac{1}{2}}, y)=0\ (i\neq0, t, -t).
\end{align*}

${\rm{(ii)}}$ We prove  $\Delta(x)=0$ for any  $x\in\mathfrak{D}.$\\[-5pt]

  For $x=\sum_{k\in\mathbb{Z}}\alpha_{k}d_{k}+\beta_{k}h_{k+\frac{1}{2}}+k_{1}{\bf c}+k_{2}{\bf l}\in\mathfrak{D}$, since $\Delta(d_{0})=\Delta(d_{1})=0$, by Lemma \ref{KM}, then $\Delta(d_{i})=0$, $i\in \mathbb{Z}$. According to Lemma \ref{3.4}, we obtain $$\Delta(x)=\lambda_{x}(\sum_{k\in\mathbb{Z}}\beta_{k}h_{k+\frac{1}{2}}+2k_{2}{\bf l})$$ and
\begin{alignat*}{2}
\Delta(x)&=&&\Delta_{d_{2t}+h_{t+\frac{1}{2}}, x}(x)
\\&=&&\big({\rm ad}(a_{2t}(d_{2t}+h_{t+\frac{1}{2}}, x)(d_{2t}+h_{t+\frac{1}{2}})+b_{0}(d_{2t}+h_{t+\frac{1}{2}}, x)h_{\frac{1}{2}}
\\[5pt]&&&+b_{-t}(d_{2t}+h_{t+\frac{1}{2}}, x)h_{-t+\frac{1}{2}}+l_{1}(d_{2t}+h_{t+\frac{1}{2}},x){\bf c}
\\&&&+l_{2}(d_{2t}+h_{t+\frac{1}{2}}, x){\bf l})+(t-\frac{1}{2})b_{-t}(d_{2t}+h_{t+\frac{1}{2}}, x)D_{1}-\frac{1}{2}b_{0}(d_{2t}+h_{t+\frac{1}{2}}, x)D_{2}\big)(x)\\&=&&a_{2t}(d_{2t}+h_{t+\frac{1}{2}}, x)(\sum_{k\in\mathbb{Z}}(2t-k)\alpha_{k}d_{2t+k}-(k+\frac{1}{2})\beta_{k}h_{2t+k+\frac{1}{2}}
\\&&&+(t+\frac{1}{2})
\alpha_{k}h_{k+t+\frac{1}{2}}+\beta_{-1-t}(t+\frac{1}{2}){\bf l}+\alpha_{-2t}\frac{4t^3-t}{6}{\bf c})\\&&&+b_{0}(d_{2t}+h_{t+\frac{1}{2}}, x)(\sum_{k\in\mathbb{Z}}\frac{1}{2}\alpha_{k}h_{k+\frac{1}{2}}+\frac{1}{2}\beta_{-1}{\bf l})
\\&&&+b_{-t}(d_{2t}+h_{t+\frac{1}{2}}, x)(\sum_{k\in\mathbb{Z}}(-t+\frac{1}{2})\alpha_{k}h_{k-t+\frac{1}{2}}+\beta_{t-1}{\bf l})
\\&&&+(t-\frac{1}{2})b_{-t}(d_{2t}+h_{t+\frac{1}{2}}, x)(\sum_{k\in\mathbb{Z}}\beta_{k}h_{k+\frac{1}{2}}+2k_{2}{\bf l})
\\&&&-\frac{1}{2}b_{0}(d_{2t}+h_{t+\frac{1}{2}}, x)(\sum_{k\in\mathbb{Z}}\alpha_{k}h_{k+\frac{1}{2}}+\beta_{-1}{\bf l})
\\&=&&a_{2t}(d_{2t}+h_{t+\frac{1}{2}}, x)(\sum_{k\in\mathbb{Z}}(2t-k)\alpha_{k}d_{2t+k}-(k+\frac{1}{2})\beta_{k}h_{2t+k+\frac{1}{2}}
\\&&&+(t+\frac{1}{2})\alpha_{k}h_{k+t+\frac{1}{2}}+\beta_{-1-t}(t+\frac{1}{2}){\bf l}+\alpha_{-2t}\frac{4t^3-t}{6}{\bf c})
\\&&&+b_{-t}(d_{2t}+h_{t+\frac{1}{2}}, x)(\sum_{k\in\mathbb{Z}}(-t+\frac{1}{2})\alpha_{k}h_{k-t+\frac{1}{2}}+\beta_{t-1}{\bf l})
\\&&&+(t-\frac{1}{2})b_{-t}(d_{2t}+h_{t+\frac{1}{2}}, x)(\sum_{k\in\mathbb{Z}}\beta_{k}h_{k+\frac{1}{2}}+2k_{2}{\bf l}).
\end{alignat*}
Comparing the two equations above, since $t\neq 0$, then we have
$$a_{2t}(d_{2t}+h_{t+\frac{1}{2}}, x)=b_{-t}(d_{2t}+h_{t+\frac{1}{2}}, x)=0.$$
Therefore  $$\lambda_{x}=0.$$
Furthermore, we have $\Delta(x)=0.$
\end{proof}

\begin{theorem}
Every 2-local derivation on the mirror Heisenberg-Virasoro algebra $\mathfrak{D}$ is a derivation.
\begin{proof}
Suppose $\Delta$ is a 2-local derivation on $\mathfrak{D}$, then there exists a derivation $\Delta_{d_{0}, d_{1}}$ such that
$$\Delta(d_{0})=\Delta_{d_{0}, d_{1}}(d_{0}), \Delta(d_{1})=\Delta_{d_{0}, d_{1}}(d_{1}).$$
Let $\Delta(1)=\Delta-\Delta_{d_{0}, d_{1}}$. Then
$$\Delta(1)(d_{0})=\Delta(1)(d_{1})=0.$$
By Lemma \ref{KM}, we have $\Delta(1)(d_{i})=0,$ $i=1,2$ and by Lemma \ref{3.4}, for fixed $t\in\mathbb{Z}\!\setminus\!\{0, 1\},$ we have
$$\Delta(1)(h_{t-\frac{1}{2}})=\lambda_{h_{t-\frac{1}{2}}}h_{t-\frac{1}{2}}.$$
Let $$\Delta(2)=\Delta(1)-\lambda_{h_{t-\frac{1}{2}}}D_{1}.$$ Then we have
\begin{align*}
&\Delta(2)(d_{0})=\Delta(1)d_{0}-\lambda_{h_{t-\frac{1}{2}}}D_{1}(d_{0})=0,\\
&\Delta(2)(d_{1})=\Delta(1)d_{1}-\lambda_{h_{t-\frac{1}{2}}}D_{1}(d_{1})=0,\\
&\Delta(2)(h_{t-\frac{1}{2}})=\Delta(1)h_{t-\frac{1}{2}}-\lambda_{h_{t-\frac{1}{2}}}D_{1}(h_{t-\frac{1}{2}})=0.
\end{align*}
According to Lemma \ref{3.5},  $\Delta(2)=0,$
 that is $$\Delta=\Delta_{d_{0}, d_{1}}+\lambda_{h_{t-\frac{1}{2}}}D_{1}.$$
Therefore, $\Delta$ is a derivation.
\end{proof}
\end{theorem}

\subsection*{Acknowledgements}
The authors are supported by the NNSF of China (Nos: 12001141, 11971134) and NSF of Hei Longjiang Province (No. JQ2020A002).

\end{document}